\documentclass[12pt]{article} 
\usepackage{amssymb} 

\begin{document} 

\title{A classification of Poisson homogeneous spaces of complex 
reductive Poisson-Lie groups} 

\author{Eugene~A.~Karolinsky\footnote{Research was 
supported in part by INTAS grant 94-4720.}} 
 
\date{\small \tt Department of Mathematics, 
Kharkov State University,\\ 
4 Svobody Sqr., Kharkov, 310077, Ukraine;\\
karol@skynet.kharkov.com; karolinsky@ilt.kharkov.ua}

\newcommand{\note}{\noindent {\bf Note.} \ }
\newcommand{\notes}{\noindent {\bf Notes.} \ } 
\newcommand{\Dg}{D(\mathfrak g)} 
\newcommand{\Ker}{\mbox{\rm Ker} \, } 
\newcommand{\Int}{\mbox{\rm Int} \, } 
\newcommand{\Aut}{\mbox{\rm Aut} \, } 
\newcommand{\Lie}{\mbox{\rm Lie} \, } 
\newcommand{\ad}{\mbox{\rm ad}} 
\newcommand{\doptr}{(\mathfrak p, \mathfrak p', \theta)}
\newcommand{\dopch}{(\mathfrak p, \mathfrak p', \theta , {\mathfrak l}_0)} 
\newcommand{\gxg}{\mathfrak g \! \times \! \mathfrak g}
\newcommand{\zxz}{\mathfrak z \! \times \! \mathfrak z'} 
\newcommand{\scpr}{\langle \cdot \, , \cdot \rangle}
\newcommand{\RRR}{{\bf R}} 
\newcommand{\UUU}{{\bf U}} 
\newcommand{\PPP}{{\bf P}} 
\newcommand{\PPPpr}{{\bf P'}} 
\newcommand{\AAA}{{\bf A}} 
\newcommand{\AAApr}{{\bf A'}} 
\newcommand{\gdiag}{{\mathfrak g}_{diag}}
\newcommand{\hdiag}{{\mathfrak h}_{diag}} 

\newtheorem{proposition}{Proposition}[section] 
\newtheorem{theorem}[proposition]{Theorem}  

\maketitle 

\begin{abstract} 

Let $G$ be a complex reductive connected algebraic group equipped 
with the Sklyanin bracket. A classification of Poisson homogeneous 
$G$-spaces with connected isotropy subgroups is given. This result is 
based on Drinfeld's correspondence between Poisson homogeneous 
$G$-spaces and Lagrangian subalgebras in the double 
$D(\mathfrak g)$ (here ${\mathfrak g} = \mbox{Lie} \, G$). 
A geometric interpretation of some of Poisson homogeneous 
$G$-spaces is also proposed.

\end{abstract} 

\begin{sloppy} 

Let $G$ be a Poisson-Lie group, $\mathfrak g = \Lie G$, 
let $\Dg$ be the double 
corresponding to the Lie bialgebra $\mathfrak g$. We say that a 
subalgebra $\mathfrak l \subset \Dg$ is called {\it Lagrangian} if it is 
a maximal isotropic subspace with respect to the natural scalar 
product in $\Dg$. It follows from \cite{Drinfeld Lagr} that there is 
a one-to-one correspondence between 
Poisson homogeneous $G$-spaces (up to isomorphism) with 
connected stabilizers and Lagrangian subalgebras 
$\mathfrak l \subset \Dg$ such that 
$\mathfrak l \cap \mathfrak g$ 
is a Lie algebra of a certain closed subgroup in $G$ 
(up to $G$-conjugacy). 

Now let $G$ be a connected complex reductive algebraic group 
equipped with the Sklyanin bracket. By $\scpr$ denote any 
nondegenerate symmetric invariant bilinear form on 
$\mathfrak g$ such that its restriction on a compact real form 
of $[\mathfrak g , \mathfrak g ]$ is positive definite. 
Then $\Dg = \gxg$, and the natural scalar product 
in $\Dg$ is given by 
   \begin{equation} \label{scpr}
     \langle (x_1, y_1), (x_2, y_2) \rangle = \langle x_1, x_2 
     \rangle - \langle y_1, y_2 \rangle ,
   \end{equation} 
where $x_1, x_2, y_1, y_2 \in \mathfrak g$ (see Section \ref{G Poisson}).  

In this paper we obtain a description of orbits of the diagonal 
$G$-action on the set of all Lagrangian subalgebras in $\gxg$ 
(see Theorem \ref{mainth}) and specify the orbits 
of Lagrangian subalgebras $\mathfrak l \subset \gxg$ such that 
the subalgebra 
$\mathfrak l \cap \gdiag \subset \gdiag \simeq \mathfrak g$  
corresponds to a certain closed subgroup in $G$ 
(see Theorem \ref{int}; here by $\gdiag$ we denote the diagonal 
image of $\mathfrak g$ in $\gxg$). 
Thus we get a classification of all Poisson homogeneous 
$G$-spaces with connected stabilizers. 

Note that the description of $G$-orbits on the set of Lagrangian 
subalgebras $\mathfrak l \subset \gxg$ such that 
$\mathfrak l \cap \gdiag = 0$ was obtained in \cite{Bel-Dr}; 
this result is related to a classification of the solutions of the 
classical Yang-Baxter equation. A classification of structures 
of a Poisson homogeneous space on $G/H$, where $H$ is 
a Cartan subgroup, was independently obtained by Jiang-Hua Lu; 
this structures are closely related to the solutions of the classical 
dynamical Yang-Baxter equation
(see \cite{Lu CDYBE}).  

This paper is organized as follows. In Section \ref{G Poisson} we recall 
the definition of the Sklyanin bracket on $G$. In Section \ref{class} 
we formulate classification theorems. Section \ref{pr} presents methods 
of the proof of Theorem \ref{mainth}. In Section \ref{geometr} we propose 
a geometric interpretation of some of Poisson homogeneous $G$-spaces, 
i.e., we construct a Poisson manifold $X$ with a Poisson $G$-action 
such that $G$-orbits on $X$ are Poisson homogeneous $G$-spaces, 
and different orbits are not isomorphic (note that in the case when  
the Poisson bracket on $G$ is zero, an analogue of $X$ is 
$\mathfrak g^*$ with the Kirillov's bracket and 
the coadjoint action of $G$). 

Note that in this paper we only formulate the main results and give 
a brief description of methods of proofs. The complete proofs 
will be presented elsewhere. 

The author is grateful to V.~G.~Drinfeld for constant attention 
to this work and to L.~L.~Vaksman and S.~Parmentier 
for useful discussions.  

\section{Poisson structure on $G$} \label{G Poisson} 

 Let us recall the definition of the Poisson structure on $G$. 
Fix a Cartan subalgebra $\mathfrak h \subset \mathfrak g$. 
Let $\RRR$ be the root system of $\mathfrak g$ with respect 
to $\mathfrak h$, $\RRR_+$ the set of positive roots 
with respect to a certain system of simple roots 
$\Gamma \subset \RRR$. Set   
$$ \mathfrak n_+ = \bigoplus_{\alpha \in \RRR_+} 
   \mathfrak g_\alpha , \ 
   \mathfrak n_- = \bigoplus_{\alpha \in \RRR_+} 
   \mathfrak g_{- \alpha}  , $$ 
$\mathfrak b_+ = \mathfrak h \oplus \mathfrak n_+$, 
$\mathfrak b_- = \mathfrak h \oplus \mathfrak n_-$. 

        Consider $r = \frac{1}{2} t_0 + t_1$ (here the tensor 
$t = t_0 + t_1 + t_2 \in \mathfrak g \otimes \mathfrak g$ 
corresponds to the bilinear form $\scpr$, 
$t_0 \in \mathfrak h \otimes \mathfrak h$, 
$t_1 \in \mathfrak n_+ \otimes \mathfrak n_-$, $t_2 \in 
\mathfrak n_- \otimes \mathfrak n_+$. 
We have $r = r_{sym} + r_{alt}$, where $r_{sym}$ is 
symmetric and $r_{alt}$ is skew-symmetric. 
Let $r^{\mu \nu}$ be the components of the tensor $r$ in some basis 
$\{ e_\mu \} \subset \mathfrak g$. Denote by $\partial_\mu$ 
(respectively by $\partial'_\mu$) the right-invariant 
(respectively left-invariant) vector field corresponding to $e_\mu$.   
Since $r$ satisfies the classical Yang-Baxter equation and 
$r_{sym}$ is $\mathfrak g$-invariant (see \cite[\S 4]{Drinfeld QG}), 
we see that the Sklyanin's formula 
$$ 
 \{ \phi , \psi \} = r^{\mu \nu} 
 (\partial'_\mu \phi \cdot \partial'_\nu \psi 
 - \partial_\mu \phi \cdot \partial_\nu \psi)   
$$ 
(here $\phi , \psi$ are regular functions on $G$)  
defines the structure of Poisson-Lie group on $G$. 
The structure of a Lie bialgebra on $\mathfrak g = \Lie G$ is defined by 
the Manin triple $(\gxg , \mathfrak g_{diag} , \mathfrak m)$, 
where $\gxg$  equipped with the scalar product (\ref{scpr}),  
$$ \mathfrak m = \{ (x,y) \in \mathfrak b_- \times \mathfrak b_+ \ | \ 
   x_\mathfrak h + y_\mathfrak h = 0 \} , $$ 
$x_\mathfrak h$ (respectively $y_\mathfrak h$) is the image of $x$ 
(respectively of $y$) in $\mathfrak h$ 
(see \cite[\S 3, Example 3.2]{Drinfeld QG}).     
In particular, the double $\Dg$ is equal to $\gxg$ equipped with  
the scalar product (\ref{scpr}).

\section{Classification theorems} \label{class} 

   Fix a Cartan subalgebra $\mathfrak h \subset \mathfrak g$.
Let ${\bf R}$ be the root system of $\mathfrak g$ with respect to 
$\mathfrak h$.

  Let ${\bf P}, {\bf P'} \subset {\bf R}$ be parabolic subsets  
(see \cite[Ch.6, \S 1.7]{Bourbaki}). Set 
$$\mathfrak p = \mathfrak h \oplus (\bigoplus_{\alpha \in {\bf P}} 
{\mathfrak g}_{\alpha}), \ 
\mathfrak p' = \mathfrak h \oplus (\bigoplus_{\alpha \in {\bf P'}} 
{\mathfrak g}_{\alpha}) .$$ 
Then $\mathfrak p$ and $\mathfrak p'$ are parabolic subalgebras 
in $\mathfrak g$. Set 
${\bf A} = {\bf P} \cap (- {\bf P}),\  
{\bf A'} = {\bf P'} \cap (- {\bf P'})$.  
Let $\mathfrak a$ and $\mathfrak a'$ be the semisimple subalgebras in 
$\mathfrak g$ generated by  
${\bf A}$ and ${\bf A'}$ respectively. Let $\tilde{\mathfrak h} = 
\mathfrak a \cap \mathfrak h ,\  
\tilde{\mathfrak h}' = \mathfrak a' \cap \mathfrak h$. 
Note that $\tilde{\mathfrak h}$ (respectively $\tilde{\mathfrak h}'$) 
is the linear span of the coroots $\alpha^\vee \in \mathfrak h$  
such that $\alpha \in \AAA$ (respectively $\alpha \in \AAApr$). 

  Let $\sigma : {\bf A} \rightarrow {\bf A'}$ be an isomorphism of 
the root systems such that $\sigma$ preserves the scalar product. 
Set 
$${\bf U} = \{ \alpha \in {\bf A'} \ | \ \sigma^{-k}(\alpha) \in
  {\bf A'} \ \ \forall \, k \in \mathbb N \} .$$ 
Since the sets $\AAA , \AAApr \subset \RRR$ are finite, and 
$\sigma : \AAA \to \AAApr$ is a bijection, we have  
$$ \UUU = \{ \alpha \in \AAA \ | \ \sigma^k (\alpha) \in \AAA \ \ 
   \forall \, k \in \mathbb N \} = $$ 
$$ = \{ \alpha \in \AAA \cap \AAApr \ | \ \sigma^l (\alpha) \in 
   \AAA \cap \AAApr \ \ \forall \, l \in \mathbb Z \} . $$ 
It is easy to prove that  
$$ 
 \mathfrak u = \mathfrak h \oplus 
 (\bigoplus_{\alpha \in {\bf U}} {\mathfrak g}_{\alpha}) 
$$ 
is a Levi subalgebra in $\mathfrak g$ (i.e., a reductive Levi subalgebra 
of a certain parabolic subalgebra in $\mathfrak g$), and ${\bf U}$ is the  
root system of $\mathfrak u$.  
We consider only the case when $\sigma$ 
preserves a certain system of simple roots in ${\bf U}$.

  Let $\xi : \mathfrak a \rightarrow \mathfrak a'$ be an isomorphism   
such that $\xi ({\mathfrak g}_{\alpha}) = {\mathfrak g}_{\sigma (\alpha)}$ 
for all $\alpha \in {\bf A}$; then $\xi (\tilde{\mathfrak h}) = 
\tilde{\mathfrak h}'$, and $\xi$ preserves $\scpr$. 
Let the linear map 
$\sigma^\vee : \tilde{\mathfrak h} \to \tilde{\mathfrak h}'$ be given by  
$\sigma^\vee (\alpha^\vee) = \sigma (\alpha)^\vee$, 
where $\alpha \in \AAA$;  
then $\xi (x) = \sigma^\vee (x)$ for all $x \in \tilde{\mathfrak h}$.  
Note that  
$$ {[\mathfrak u , \mathfrak u]}^{\xi} = 
   \{ x \in [\mathfrak u , \mathfrak u] \ | \ 
   \xi (x) = x \} $$ 
is a reductive Lie algebra, and ${\mathfrak h}^{\xi} = 
{[\mathfrak u , \mathfrak u]}^{\xi} \cap \mathfrak h$ is  
a Cartan subalgebra in ${[\mathfrak u , \mathfrak u]}^{\xi}$ 
(see \cite[Ch.4, \S 4.2]{V-O}). 

  Consider a nilpotent element 
$x \in {[\mathfrak u , \mathfrak u]}^{\xi}$  
(we say that   
an element $x \in \mathfrak g$ is called nilpotent if 
$x \in [\mathfrak g , \mathfrak g]$ and $\mbox{ad} x$ is nilpotent). 
Let $h \in {\mathfrak h}^{\xi}$ be the characteristic   
of the nilpotent element $x$ (see   
\cite[Ch.6, \S 2.1]{Lie 3}; recall that one can reconstruct $x$ by $h$ 
uniquely up to conjugation). Let the isomorphism 
$\theta : \mathfrak a \rightarrow \mathfrak a'$ be given by 
$\theta = \xi \cdot \exp(\mbox{ad} x)$.

  Let $\mathfrak z$ (respectively $\mathfrak z'$) be the orthogonal 
complement to $\tilde{\mathfrak h}$ 
(respectively $\tilde{\mathfrak h}'$) 
in $\mathfrak h$. Note that the natural maps  
$\mathfrak z \to \mathfrak p / [\mathfrak p , \mathfrak p]$ and  
$\mathfrak z' \to \mathfrak p' / [\mathfrak p' , \mathfrak p']$ 
are isomorphisms. 
Consider  
$\zxz$ equipped with the scalar product (\ref{scpr}). 
Let ${\mathfrak l}_{0} \subset \zxz$ be a Lagrangian subspace. 
Consider  
$$ 
  \mathfrak l = \{ (x,y) \in \gxg \ | \ x \in \mathfrak p , 
  \ y \in \mathfrak p' , \ 
  \theta (x_{\mathfrak a}) = y_{\mathfrak a'} , \ 
  (x_{\mathfrak z} , y_{\mathfrak z'}) \in 
  {\mathfrak l}_{0} \} ;     
$$ 
here $x_{\mathfrak a}$ is the image of $x$ in $\mathfrak a$, \ 
$x_{\mathfrak z}$ is the image of $x$ 
in $\mathfrak z = \mathfrak p / [\mathfrak p , \mathfrak p]$, \ 
$y_{\mathfrak a'}$ is the image of $y$ 
in $\mathfrak a'$, \ $y_{\mathfrak z'}$ is the image of $y$ in  
$\mathfrak z' = \mathfrak p' / [\mathfrak p' , \mathfrak p']$. 
Then $\mathfrak l$ is a Lagrangian subalgebra in $\gxg$.  
By $L({\bf P}, {\bf P'}, \sigma , \xi , h, {\mathfrak l}_{0})$ denote 
the class of $G$-conjugacy of $\mathfrak l$. 
  
\begin{theorem} \label{mainth}
  {\rm 1)} Any $G$-orbit on the set of all 
  Lagrangian subalgebras in $\gxg$ is of the form 
  $L({\bf P}, {\bf P'}, \sigma , \xi , h, {\mathfrak l}_{0})$.
  
  {\rm 2)} $L({\bf P}, {\bf P'}, \sigma , \xi , h, {\mathfrak l}_{0}) = 
  L(\tilde{\bf P}, \tilde{\bf P'}, \tilde{\sigma}, \tilde{\xi}, 
  \tilde{h}, \tilde{\mathfrak l}_{0})$ iff  
  $({\bf P}, {\bf P'}, \xi , h, {\mathfrak l}_{0})$ and $(\tilde{\bf P}, 
  \tilde{\bf P'}, \tilde{\xi}, \tilde{h}, \tilde{\mathfrak l}_{0})$ are 
  $N(\mathfrak h)$-conjugate (here by $N(\mathfrak h)$ we denote 
  the normalizer of $\mathfrak h$ in $G$). 
\end{theorem} 

\notes 
\rm 1) Let $W$ be the Weyl group of the root system $\RRR$. 
If $({\bf P}, {\bf P'}, \xi , h, 
{\mathfrak l}_{0})$ and  
$(\tilde{\bf P}, \tilde{\bf P'}, \tilde{\xi}, \tilde{h},
\tilde{\mathfrak l}_{0})$ are $N(\mathfrak h)$-conjugate, 
then $({\bf P}, {\bf P'}, \sigma , h)$ and  
$(\tilde{\bf P}, \tilde{\bf P'}, \tilde{\sigma}, \tilde{h})$ 
are $W$-conjugate. 

2) Every class of $G$-conjugacy 
of Lagrangian subalgebras in $\gxg$ 
depends on the discrete parameters $({\bf P}, {\bf P'}, \sigma , h)$ 
and the continuous parameters $(\xi , {\mathfrak l}_{0})$. Fix  
$({\bf P}, {\bf P'}, \sigma)$ and denote by $\Xi$ (respectively by  
$\Lambda$) the space of parameters $\xi$ 
(respectively ${\mathfrak l}_{0}$)  
such that $\xi$ (respectively ${\mathfrak l}_{0}$) corresponds to 
$({\bf P}, {\bf P'}, \sigma)$. 
Let $\langle {\bf A} \rangle$ be the linear span of ${\bf A}$,   
$n = \dim \mathfrak z$.  
It can be proved that  
$$\dim \Xi = \dim \ \{ \alpha \in \langle {\bf A} \rangle \ | \ 
  \sigma (\alpha ) = \alpha \} , $$ 
$\dim \Lambda = \frac{n(n-1)}{2}$ 
\ (note that $\Lambda$ is the Lagrangian 
Grassmann manifold for $\zxz$). 

\bigskip  

   We shall say that a class of $G$-conjugacy  
$L (\PPP , \PPPpr , \sigma , \xi , h , \mathfrak l_0)$ is called 
{\it integrable} 
(respectively {\it algebraic integrable})  
if the subalgebra 
$\mathfrak l \cap \mathfrak g_{diag} \subset \mathfrak g_{diag} 
\simeq \mathfrak g$ corresponds to a closed (respectively  
closed by Zariski) subgroup in $G$ 
for a certain (and then for every) Lagrangian subalgebra 
$\mathfrak l \in L (\PPP , \PPPpr , \sigma , \xi , h , \mathfrak l_0)$. 
Theorem \ref{int} gives a test of the integrability  
and the algebraic integrability of  
$L (\PPP , \PPPpr , \sigma , \xi , h , \mathfrak l_0)$. 

    Let $H \subset G$ be the connected subgroup such that 
$\Lie H = \mathfrak h \subset \mathfrak g$.    

\begin{theorem} \label{int} 
  A class of $G$-conjugacy 
  $L (\PPP , \PPPpr , \sigma , \xi , h , \mathfrak l_0)$ 
  is integrable (respectively algebraic integrable) iff  
  the subspace 
  $$ V = \{ x \in \mathfrak h \ | \ 
     (x_\mathfrak z , x_{\mathfrak z'}) \in \mathfrak l_0 , \  
     \sigma^\vee (x_{\tilde{\mathfrak h}}) = x_{\tilde{\mathfrak h}'} \} 
     \subset \mathfrak h $$ 
  (here $x_\mathfrak z$ is the image of $x$ in $\mathfrak z$, ets.) 
  is the Lie algebra of a closed   
  (respectively closed by Zariski) subgroup in $H$.  
\end{theorem}  

\note 
\rm It follows from the Theorem \ref{int} that the (algebraic) 
integrability of a class of $G$-conjugasy  
$L (\PPP , \PPPpr , \sigma , \xi , h , \mathfrak l_0)$  
depends only on $\sigma$ and $\mathfrak l_0$ 
(and is independent of $\xi$ and $h$). 

\bigskip 

  Now we recall a well-known method to verify that a subspace  
$V \subset \mathfrak h$ is the Lie algebra of a closed 
(respectively closed by Zariski) subgroup in $H$. 

  Consider the lattice 
$$ {\cal H} = \Ker (\exp : \mathfrak h \to H) \subset \mathfrak h. $$  

\begin{proposition} \label{V alg int} 
  {\rm (see {\cite[Ch.3, \S 2, Theorem 5]{V-O}})} 
  A subspace $V \subset \mathfrak h$ corresponds to a closed 
  by Zariski subgroup in $H$ iff  
  $V$ is defined over $\mathbb Q$ with respect 
  to the lattice $\cal H$, 
  i.e., $V = \cal V \otimes \mathbb C$ for a certain sublattice  
  $\cal V \subset \cal H$.  
\end{proposition} 

Let $\mathfrak t = \cal H \otimes \mathbb R \subset \mathfrak h$.  

\begin{proposition} \label{V int} 
  A subspace $V \subset \mathfrak h$ corresponds to a 
  closed subgroup in $H$ iff $V \cap \mathfrak t$ 
  is defined over $\mathbb Q$ with respect to the lattice $\cal H$, 
  i.e., $V \cap \mathfrak t = \cal V \otimes \mathbb R$ for a certain  
  sublattice $\cal V \subset \cal H$. 
\end{proposition}

\section{Methods of the proof of Theorem \ref{mainth}} \label{pr} 

Now we present a way to proove Theorem \ref{mainth}. 

Let $\mathfrak p, \mathfrak p' \subset \mathfrak g$ 
be parabolic subalgebras. We have 
$\mathfrak p / \mathfrak p^\perp = 
\mathfrak a \oplus \mathfrak z$, 
where $\mathfrak a$ is semisimple, 
and $\mathfrak z$ is abelian; 
the same holds for $\mathfrak p'$.  
Let $\theta : \mathfrak a \rightarrow \mathfrak a'$ 
be an isomorphism such that $\theta$ preserves $\scpr$.  
We shall say that a triple $\doptr$ is called {\it admissible}. 
By $T(\mathfrak g)$ denote the set of all admissible triples.  

Consider $\doptr \in T(\mathfrak g)$. Let 
${\mathfrak l}_0 \subset \mathfrak z \! \times \! \mathfrak z'$ 
be a Lagrangian subspace with respect to the bilinear form 
(\ref{scpr}). We say that a quadruple $(\mathfrak p , \mathfrak p',
\theta , {\mathfrak l}_0)$ is called {\it admissible}. Suppose 
$(\mathfrak p , \mathfrak p' , \theta , {\mathfrak l}_0)$ 
is an admissible quadruple; then set  
$$\mathfrak l (\mathfrak p , \mathfrak p' , 
  \theta , {\mathfrak l}_0) := \{ (x , y)
  \in \mathfrak p \! \times \! \mathfrak p' \  | \ 
  \theta ( x_{\mathfrak a} ) = y_{\mathfrak a'} ,\ 
  ( x_{\mathfrak z} , y_{\mathfrak z'} ) \in {\mathfrak l}_0 \} \subset 
  \mathfrak g \! \times \! \mathfrak g ,$$ 
where $x_{\mathfrak a}$ is the image of $x$ in $\mathfrak a$,
$x_{\mathfrak z}$ is the image of $x$ in $\mathfrak z$, 
$y_{\mathfrak a'}$ 
is the image of $y$ in $\mathfrak a'$, $y_{\mathfrak z'}$ 
is the image of $y$ in $\mathfrak z'$. 
It is not hard to prove the following proposition. 

\begin{proposition} \label{dopch1}  
  {\rm 1)} $\mathfrak l (\mathfrak p , \mathfrak p' , \theta , 
     {\mathfrak l}_0)$ is a Lagrangian subalgebra. 

  {\rm 2)} The correspondence $\dopch \mapsto \mathfrak l \dopch$ 
     is a $G$-equivariant bijection between the set of all Lagrangian 
     subalgebras in $\gxg$ and the set of all admissible quadruples 
     $(\mathfrak p , \mathfrak p' , \theta , {\mathfrak l}_0)$. 

  {\rm 3)} Lagrangian subalgebras $\mathfrak l (\mathfrak p , 
     \mathfrak p' , \theta , {\mathfrak l}_0)$ and 
     $\mathfrak l (\mathfrak p , \mathfrak p' , \theta , 
     {\tilde{\mathfrak l}}_0)$ are $G$-conjugate iff 
     ${\mathfrak l}_0 = {\tilde{\mathfrak l}}_0$.
\end{proposition}

Thus a classification of Lagrangian subalgebras is reduced to 
a classification of admissible triples up to $G$-conjugacy. 
It can be shown 
that the theory of admissible triples is quite similar to the theory of 
automorphisms of complex semisimple Lie algebras. 
In fact, there exists 
a natural notion of a semisimple admissible triple; 
we can define a notion 
of an invariant subalgebra for an admissible triple; 
for any semisimple admissible triple there exists an invariant Cartan 
subalgebra; it is possible, using invariant Cartan subalgebras, to give 
a complete description of semisimple admissible triples up to 
$G$-conjugacy; for any admissible triple there exists an analogue of 
the Jordan decomposition, etc. The realization of this program 
leads us to Theorem \ref{mainth}.

\section{A geometric interpretation} \label{geometr} 

In this section we give a geometric interpretation of some 
of Poisson homogeneous $G$-spaces. 

   By $\bar{G}$ denote the group of all automorphisms 
$g : \mathfrak g \to \mathfrak g$ such that 
the following conditions hold: 
(1) $g$ preserves the scalar product $\scpr$;  
(2) $g$ is equal to the identity mapping on the 
center of $\mathfrak g$. 
Suppose $g \in \bar{G}$ and set  
$$ \mathfrak l_g = \{ (x,y) \ | \ x = g(y) \} \subset \gxg . $$ 
Then $\mathfrak l_g$ is a Lagrangian subalgebra. 
Note that the Lagrangian 
subalgebras $\mathfrak l_g$ form the classes of $G$-conjugacy  
$L(\PPP , \PPPpr , \sigma , \xi , h , \mathfrak l_0)$ such that 
$\PPP = \PPPpr = \RRR$ and $\mathfrak l_0$ is the image 
of the center of 
$\mathfrak g$ under the diagonal mapping $\mathfrak g \to \gxg$.   

        Let us give a geometric interpretation of Poisson homogeneous 
$G$-spaces corresponding to Lagrangian subalgebras 
of the form $\mathfrak l_g$. 
Note that the connected component 
of the center of $G$ acts trivially on 
the subalgebras $\mathfrak l_g$; 
therefore it is enough to consider the case 
when $G$ is semisimple, i.e., 
$[\mathfrak g , \mathfrak g] = \mathfrak g$. 
In the following part of this section we consider the case 
$G = \Int \mathfrak g$. 

Let $\phi , \psi$ be regular functions on $\bar{G}$. Consider 
\begin{eqnarray} \label{X}  
  \{ \phi , \psi \} & = & - r^{\mu \nu}_{alt} \cdot (\partial'_\mu \phi - 
  \partial_\mu \phi) \cdot (\partial'_\nu \psi - \partial_\nu \psi) 
  \nonumber \\  
  & + & r^{\mu \nu}_{sym} \cdot (\partial'_\mu \phi - 
  \partial_\mu \phi) \cdot (\partial'_\nu \psi + \partial_\nu \psi) ,     
\end{eqnarray} 
where $r$, $\partial_\mu$ and $\partial'_\mu$ are defined in Section 
\ref{G Poisson}. By $X$ we denote the manifold $\bar{G}$ 
equipped with the bracket (\ref{X}).   

\begin{theorem} \label{geom}
    The bracket {\rm (\ref{X})} is a Poisson bracket,
the action of $G$ on $X$ by conjugations is Poisson, and the orbits of 
this action are Poisson homogeneous $G$-spaces such that the 
Lagrangian subalgebra $\mathfrak l_g$ corresponds to a point 
$g \in X$.  
\end{theorem} 

\note  
\rm The bracket (\ref{X}) is a special case of the bracket from 
\cite[Theorem~3.1]{Parm}, when $J_1 = - J_2$ (using 
the notation from \cite{Parm}). See also \cite{Li Parm}. 

\bigskip 

Theorem \ref{geom} can be proved by using the following 
general result (see Theorem \ref{D/G}). 
Suppose $G$ is an arbitrary Poisson-Lie group.  
We say that {\it a double} of $G$ 
is a Lie group $D$ such that the following conditions hold: 
(1) $\Lie D = \Dg$; 
(2) The natural scalar product in $\Dg$ is invariant with respect to  
the adjoint action of $D$ (then $D$ becomes a Poisson-Lie group 
by means of the canonical element  
$r \in \mathfrak g \otimes \mathfrak g^* \subset \Dg \otimes \Dg$, 
see \cite[\S 13]{Drinfeld QG}); 
(3) $G$ is a closed Poisson-Lie subgroup in $D$. 

\begin{theorem} \label{D/G}  
 {\rm \cite{Kar2}}  
 Let $G$ be a Poisson-Lie group,  
 $\mathfrak g = {\rm \Lie} G$. Let $D$ be a double of $G$. 
 Consider the action of $G$ on the Poisson manifold $D/G$ by 
 left translations.
 Suppose $w \in D$ and denote by $x$ the image of $w$ in $D/G$; 
 then $X = G \cdot x$ is a Poisson homogeneous $G$-space, and 
 the Lagrangian subalgebra  
 $\mathfrak l_x := w \cdot \mathfrak g \cdot w^{-1} \subset \Dg$ 
 corresponds to the pair $(X, x)$. 
\end{theorem}

In our case take $D = \bar{G} \times \bar{G}$; 
then Theorem \ref{geom} 
follows from Theorem \ref{D/G}.  

\end{sloppy}

\end{document}